\documentclass[12pt,twoside,a4paper]{amsart}
\usepackage{amssymb,url,verbatim,hyperref,wasysym}
\date{\today}

\newenvironment{nouppercase}{%
  \renewcommand{\uppercasenonmath}[1]{}}{}

\def \w {\wedge}
\def \contr {\mathrel{\lrcorner}}

\def \C {\mathbb{C}}

\def \ideala {\mathfrak{a}}
\def \dbar {\overline{\partial}}
\def \O {\mathcal{O}}
\def \local {\O_{\C^n,0}}

\def \localz {\O_{Z,0}}
\def \localx {\O_{X,0}}
\def \localv {\O_{V,0}}

\def \idealJ {\mathcal{J}}

\def \pbpi {{\pi'}^*}
\def \pb {\pi^*}
\def \dmod {\mathcal{N}(Z)}
\def \sep {\text{SEP}}

\def \chx {\mathcal{CH}_X}
\def \vec {\mathcal{S}}

\DeclareMathOperator{\End}{End}

\DeclareMathOperator{\Image}{Im}
\DeclareMathOperator{\Hom}{Hom}
\DeclareMathOperator{\ann}{ann}

\DeclareMathOperator{\order}{ord}
\DeclareMathOperator{\rank}{rank}
\DeclareMathOperator{\kersh}{\mathcal{K}er}

\def\ba{\begin{align*}}
\def\ea{\end{align*}}
\def\bal{\begin{align}}
\def\eal{\end{align}}

\newtheorem{theorem}{Theorem}[section]
\newtheorem{lemma}[theorem]{Lemma}

\newtheorem{proposition}[theorem]{Proposition}
\theoremstyle{definition}
\newtheorem{example}[theorem]{Example}
\newtheorem{definition}[theorem]{Definition}
\theoremstyle{remark}
\newtheorem{remark}[theorem]{Remark}

\title[]{A Brian\c con-Skoda type result for a non-reduced analytic space}

\begin{document}

\author{Jacob Sznajdman}

\address{Jacob Sznajdman\\Mathematical Sciences, Chalmers University and Gothenburg University\\S-412 96 GOTHENBURG\\SWEDEN}

\email{sznajdma@chalmers.se}

\subjclass[2000]{32C30, 13B22, 32A26}

\begin{abstract}
  We present here an analogue of the Brian\c con-Skoda theorem for a germ
  of an analytic space $Z$ at $0$, such that $\localz$ is not necessarily reduced.
\end{abstract}

\begin{nouppercase}
\maketitle
\end{nouppercase}

\bibliographystyle{amsplain}

\section{Introduction}\label{intro}
The Brian\c con-Skoda theorem, \cite{brianconskoda}, states that for any
ideal $\ideala \subset \local$ generated by $m$ germs, we have the inclusion
$\overline{\ideala^{\min(m,n)+r-1}}\subset\ideala^r$, where $\overline{I}$
denotes the integral closure of $I$. We will refer to this theorem as the
classical Brian\c con-Skoda theorem.
The generalization to an arbitrary regular Noetherian ring was proven
algebraically in \cite{lipmansathaye}.

Huneke, \cite{huneke:uniform}, showed that for a quite general Noetherian reduced local ring $S$
there is an integer $N$ such that $\overline{\ideala^{N+r-1}}\subset\ideala^r$
for all ideals $\ideala \subset S$ and $r\geq 1$. In particular this applies when $S=\localv$,
the local ring of holomorphic functions of a germ of a reduced analytic space $V$.
This case of Huneke's theorem was recently reproven analytically, \cite{ass:bsvar}.
Assume that $\ideala$ is generated by some elements $a_i$ and let $|\ideala|^2 = \sum_1^m |a_i|^2$; up to constants
this does not depend on the choice of the generators.
Since a function $\phi$ in
$\overline{\ideala^M}$ is characterized by the property that $|\phi|\leq C|\ideala|^M$,
\cite{lejeunejalaberttessier}, an equivalent formulation of the theorem is that
$\phi$ belongs to $\ideala^r$ whenever $|\phi|\leq C|\ideala|^{N+r-1}$ (on $V$).

We will consider a germ of an analytic space, that is, a pair $Z=(X,\localz)$
of a germ of a reduced analytic variety $X \subset \C^n$ at $0\in \C^n$ and its local ring $\localz = \local / \idealJ$,
where $\idealJ \subset \local$ is an ideal such that $Z(\idealJ) = X$. 
We assume throughout this paper that $\idealJ$ has pure dimension. 

The aim of this paper is to find an appropriate generalization of the Brian\c con-Skoda
theorem to this setting -- when $S=\localz$.
Now that we have dropped the assumption that $Z$ is reduced, the situation becomes different;
the integral closure of any ideal contains the nilradical $\sqrt{0}$ by definition, so 
$\overline{\ideala^N}\subset\ideala$ can only hold if
$\sqrt{0} \subset \ideala$.
In the following example we consider the most simple non-reduced
space, which will help to illustrate some general notions and our main result Theorem~\ref{bs}.

\begin{example}\label{enkeltex}
  Consider the analytic space $Z$ such that 
  \begin{align*}
    X&=\left\{w=0\right\}=\C^{n-1}\subset\C^n \\
    \O_Z &= \C[[z_1,\dots z_{n-1},w]]/(w^k),\quad k\geq 2.
  \end{align*}
  The nilradical is $(w)$, and is not contained in $\ideala=(w^2)$ if $k>2$. It may be helpful to think of
  the space $Z$ as $\C^{n-1}$ with an extra infinitesimal direction transversal to $X$,
  and its structure sheaf being the $k$:th order Taylor expansions in that direction. For each $f\in\localz$ we have
  \begin{align*}
    f(z,w) &= \sum_{i=0}^{k-1} \frac{\partial^i f}{\partial w^i}(z,0) \frac{w^i}{i!}, \quad \text{and} \\
    \localz &\simeq 
                     \O_{X,0}^{\oplus k}.
  \end{align*}
  Although the function $w$ is identically zero on $X$ (so that $|w| \leq C|\ideala|^M$ on $X$ for any $M$),
  the element $w$ does not belong to $\ideala = (w^2)$.
  Since $Z$ is non-reduced, evaluating $w$, or any other element, as a function on $X$ does not give enough
  information to determine ideal membership. We also have to take into account the transversal derivatives. 
\end{example}

A germ of a holomorphic differential operator $L$ is called
\textit{Noetherian} with respect to an ideal $\idealJ \subset \local$
if $L\phi \in \sqrt{\idealJ}$ for all $\phi \in \idealJ$.
We say that $L_1, \dots ,L_M$ is a \textit{defining} set of Noetherian operators for $\idealJ$
if $\phi \in \idealJ$ if and only if $L_1 \phi, \dots, L_M \phi \in \sqrt{\idealJ}$.
The existence of a defining set for any ideal $\idealJ$ is due to Ehrenpreis \cite{ehrenpreis}
and Palamodov \cite{palamodov}, see also
\cite{bjork:abel}, \cite{hormander:scv} and \cite{oberst}. As for the example above,
$1, \partial/\partial w, \dots, \partial^{k-1}/\partial w^{k-1}$ is a defining set for $(w^k)$.

If $L$ is Noetherian with respect to $\idealJ$, then $L \psi$ is a well-defined function on
$X$ for any $\psi \in \localz$, and $L$ induces a mapping
$L: \localz \to \O_{X,0}$. Let $\dmod$ be the set of all such mappings;
this set does not depend on the local representation of $Z$ as a subscheme of $\C^n$.
If $(L_i)$ is a defining set for $\idealJ$, then by definition any element $\phi\in\localz$ is determined
uniquely by the tuple of functions $(L_i \phi)$ on $X$, cf. Example~\ref{enkeltex}.
This fact indicates that it is natural to impose size conditions on the whole set $(L_i \phi)$ to generalize the
Brian\c con-Skoda theorem:

\begin{theorem}\label{bs}
  Let $Z$ be a germ of an analytic space 
  such that $\mathcal{O}_Z$ has pure dimension.
  Then there exists an integer $N$ and operators $L_1,\dots,L_M \in \dmod$ such that
  for all ideals $\ideala\subset \mathcal{O}_Z$ and all $r\geq 1$,
  \begin{align} \label{premise}
    |L_j \phi|\leq C |\ideala|^{N+r-1} \,\,\, \textrm{on} \,\,\, X, \quad 1 \leq j \leq M, 
  \end{align}
  implies that $\phi\in\ideala^r$.
\end{theorem}

In Section \ref{smooth} we give a version of Theorem~\ref{bs} for the case when $X$ is smooth.

Andersson and Wulcan gave in \cite{aw:semester} a proof of a global version of the Brian\c con-Skoda-Huneke theorem
on a reduced singular variety, which is a version of the effective Nullstellensatz. The proof is
based on the corresponding local result in \cite{ass:bsvar}. One can therefore hope that
Theorem~\ref{bs} can be globalized in a similar way, cf., Remark~\ref{global} below.

Although the formulation of Theorem \ref{bs} is intrinsic, we will choose an embedding and
work in the ambient space exclusively. 
For the remainder of this paper we fix a choice of functions $a_j \in \local$, $1 \leq j \leq m$,
 so that the images of $a_j$
in $\localz$ generate $\ideala$. We will also identify $\phi \in \localz$ with an arbitrary representative in $\local$,
and each operator $L \in \dmod$ with an operator on $\C^n$ that represents it.
With this point of view, we are to show that if the inequalities \eqref{premise} hold on $X \subset \C^n$,
then (the representative of) $\phi$ belongs to $(a)^r + \idealJ$, where $(a)=(a_1,\ldots,a_m)$.


\begin{example}\label{bevisatexempel}
  In this example we shall give a direct proof of Theorem~\ref{bs}, based on the classical Brian\c con-Skoda theorem,
  for the analytic space in Example~\ref{enkeltex}. For simplicity, we will assume that $r=1$.
  As we saw above,
  a set of defining differential operators for the ideal $(w^k)$ is formed by $L_j = \partial^j/\partial w^j$,
  $0\leq~j~\leq k-1$. Let
  \begin{align*}
    \varrho&=\min(n-1,m),\ \quad \text{and}\\
    |a|^2 &= \sum_1^m |a_i|^2.
  \end{align*}
  We will prove that if
  \begin{align}\label{smpremise}
    |L_j \phi| \lesssim |a|^{\varrho + k-1-j},
  \end{align}
  on $X = \{w=0\}$, then $\phi \in (a) + (w^k)$. Thus Theorem~\ref{bs} follows in this case with
  $N = \varrho + k - 1$. We will also show that none of the hypotheses \eqref{smpremise} can be relaxed.

  In the proof, we will allow ourselves to abuse notation; for example, we will write simply $(a)$ when we
  actually are referring to some element that belongs to $(a)$.
  From \eqref{smpremise} and the classical Brian\c con-Skoda theorem 
for ${\O_{X,0}} = {\O_{\C^{n},0} / (w)}$, we get
  \begin{align}\label{phipartials}
    \frac{\partial^j\phi}{\partial w^j} = (a)^{k-j} + k_j w, \quad k_j \in \local.
  \end{align}
  We will show inductively that
  \begin{align} \label{induktion}
    \phi = \sum_{i=0}^p w^i (a)^{k-i} + g_p w^{p+1}, \quad g_p \in \local,
  \end{align}
  holds for $p \leq k-1$.
  First assume that $p=0$. Then
  \eqref{induktion} reduces to
  \eqref{phipartials} with $j=0$. Now assume that \eqref{induktion} holds for some $p < k-1$.
  Let us differentiate \eqref{induktion} $p+1$ times with respect to $w$, and compare the result with
  \eqref{phipartials} for $j=p+1$.
  This gives
  \begin{align*}
    g_p \in (a)^{k-p-1} + (w).
  \end{align*}
  If we substitute this back into \eqref{induktion}, we get
  \begin{align*}
    \phi \in \sum_{i=0}^{p+1}w^i (a)^{k-i} + (w^{p+2}).
  \end{align*}
  This completes the proof of \eqref{induktion}.
  Now note that for $p=k-1$, \eqref{induktion} implies that $\phi \in (a) + (w^k)$, which completes the proof.

  Finally, we will show that if any of the hypotheses \eqref{smpremise} are relaxed, it is possible that $\phi \notin (a) + (w^k)$.
  To this end, we need to find $\phi_p$ for each $0\leq p \leq k-1$,
  such that $\phi_p \notin (a)+(w^k)$ and
  \begin{align}\label{phi_p}
    |\partial^j_w \phi_p| &\leq |a|^{\varrho +k-1-j}, \quad j\neq p \\\label{phi_p_2}
    |\partial^p_w \phi_p| &\leq |a|^{\varrho +k-2-p}.
  \end{align}
  Take $n=2$ and $(a)=(z+w)$. Then $\varrho=1$, since $n-1=m=1$.
  A suitable choice is now $\phi_p = w^p z^{k-1-p}$. It is easy to verify \eqref{phi_p}-\eqref{phi_p_2}.
  The function $\phi_p$ does not belong to $(a)+(w^k)$, because if it did
  we would have
  \begin{align*}
    w^p z^{k-1-p} - (z+w)(a_0(z) + \dots + a_{k-1}(z)w^{k-1}) \in (w^k),
  \end{align*}
  which would give
  $a_0=\dots=a_{p-1}=0$ and $za_p = z^{k-1-p}$, $za_{p+1}=-a_p$, $za_{p+2}=-a_{p+1}$, etc, so
  $a_{k-1}=\pm 1/z$. This is a contradiction since $a_{k-1}$ is holomorphic at $0\in \C^n$.
\end{example}

The idea of the proof of Theorem~\ref{bs} is to use a certain residue current associated to $(a)^r$ and $\idealJ$; if $\phi\in\local$ annihilates this
current, then by solving a sequence of $\dbar$-equations it follows that $\phi$
belongs to $(a)^r+\idealJ$, so that the image of $\phi$ in $\localz$ belongs to $\ideala^r$.
Alternatively, one
can also use a division formula to obtain an explicit integral representation of the membership $\phi \in (a)^r + \idealJ$
whenever $\phi$ annihilates the associated residue current.
These ways of proving ideal membership are used in \cite{ass:bsvar} and go back to \cite{andersson:residueideals}
and \cite{andersson:strong}.
We will show in Section~\ref{proof} that $\phi$ annihilates the residue current mentioned above whenever
\eqref{premise} holds.

\section{Coleff-Herrera currents and Noetherian differential operators}\label{coleffherrera}
Assume that $X$ is a germ of an analytic set of pure codimension $p$ at $0\in\C^n$.
Let $\mu$ be a current of bidegree $(0,p)$ with support on $X$.
Throughout this paper we let $\chi$ be a smooth function such that $\chi \equiv 0$ on $[0,x_1]$ and $\chi \equiv 1$ on $[x_2,\infty)$
for some $0<x_1 < x_2<\infty$.
  One says that $\mu$ has the standard extension property ($\sep$)
  if $\mu~ =~\lim_{\varepsilon \to 0} \chi(|h|^2/\varepsilon^2)\mu$ for any
  $h \in \local$ that does not vanish identically on any component of $X$.
  \begin{definition}
    A current of bidegree $(0,p)$ with support on $X$ is a Coleff-Herrera current on $X$ if
    it is $\dbar$-closed, has the $\sep$ and is annihilated by $\overline{\phi}$ 
for all  $\phi \in \local$ that vanish on $X$. 
  \end{definition}
  The set of all Coleff-Herrera currents on $X$ is an $\local$-module which we denote by $\chx$.
  For any $\mu \in \chx$, $\ann \mu$ is a pure-dimensional ideal whose associated
  primes correspond to the irreducible components of $X$. 
  In Theorem~\ref{bjork_thm} below, which is due to Bj\"ork, \cite{bjork:abel}, we have principal value integrals of the form
  \begin{align}\label{pv_def}
    \int_X \frac{\xi}{h^{k}} := \lim_{\varepsilon \to 0} \int_X \frac{\chi(|h|^2/\varepsilon^2)\xi}{h^{k}},
  \end{align}
  where $h\in \localx$ and $\xi$ is a test form. For details about the existence of such integrals,
  see e.g. \cite{bjorksamuelsson:reg}. The right hand side of \eqref{pv_def} is independent of the choice of
  $\chi$ as long as it has the properties mentioned above.
  We will use the symbol $\contr$ for interior multiplication of a differential form
  by a section of $\bigwedge^p T(\C^n)$.

  \begin{theorem}[Bj\"ork, \cite{bjork:abel}]\label{bjork_thm}
    Given $\mu \in \chx$, there is a multi-index $M$, two finite sets of holomorphic differential operators
    $\{L_\alpha\}$ and $\{K_\alpha\}$ for $\alpha \leq M$, an integer $N_0$, an analytic function
    $h$, and a smooth section $\vec$ of $\bigwedge^p T(\C^n)$, such that $\{L_\alpha\}$ is a defining set for the ideal
    $\ann \mu$ and
    \begin{align}\label{phimu'}
      \phi\mu.\xi = \int_X \sum_{\alpha \leq M} \frac {1}{h^{N_0}} L_\alpha (\phi) K_\alpha (\vec \contr \xi),
    \end{align}
    for any holomorphic function $\phi$ and test form $\xi$.
    \begin{proof}
      The proof is essentially that of Bj\"{o}rk, \cite{bjork:abel}, but
      we include it for the reader's convenience since our formulation of the theorem is slightly different.
      It follows from the local parametrization theorem that one can find holomorphic functions
      $f_1, \dots, f_p$ forming a complete intersection, such that $X$ is a union of a number of
      irreducible components of $V_f = \{f_1=\dots=f_p=0\}$ and $df_1 \wedge \ldots \wedge df_p \neq 0$
      generically on $X$. Let the first $n-p$ coordinates in $\C^n$ be denoted by
      $\zeta$ and the last $p$ ones by $\eta$.
      We then get local holomorphic coordinates
      \begin{align}
        z &=\zeta \notag\\
        w &=f(\zeta,\eta)\label{coordinates}
      \end{align}
      outside the hypersurface $W$ defined by
      \begin{align}\label{detdef}
        h:=\det \frac{\partial f}{\partial \eta}.
      \end{align}
      By possibly rotating the coordinates $(\zeta,\eta)$, we can make sure that $h$ will not vanish identically on
      any irreducible component of $X$. 
      
      We will first show that
      \eqref{phimu'} holds for some operators $L_\alpha$ when $\xi$ has support outside of $\{h=0\}$.
      We will then see that it follows from this special case of \eqref{phimu'} that the operators $L_\alpha$
      indeed form a defining set for $\ann \mu$. Finally, we will show that \eqref{phimu'} holds for general $\xi$.

      Since $\mu$ is a Coleff-Herrera current on the complete intersection $V_f$, we get by 
by Theorem 4.2 in \cite{andersson:ch} that
      \begin{align}\label{chproduct}
        \mu = A\left[\dbar \frac 1{f_1^{1+M_1}} \wedge \dots \wedge \dbar \frac 1{f_p^{1+M_p}} \right]
      \end{align}
      for some  integers $M_j$ and  holomorphic function $A$.
      A basic fact is that for an $(n,n-p)$ test form $\xi$,
      \begin{align*}
        \left[\dbar \frac 1{w_1^{1+M_1}} \wedge \dots \wedge \dbar \frac 1{w_p^{1+M_p}}\right].\xi
        = M! (2\pi i)^p \int_{w=0} \partial_{w_1}^{M_1} \dots \partial_{w_p}^{M_p} \left(\frac{\partial}{\partial w} \contr \xi\right),
      \end{align*}
      where the derivative symbols refer to Lie derivatives and $\partial/\partial w = \partial/\partial w_1 \wedge \ldots \wedge \partial/\partial w_p$.
      Using first Leibniz' rule and \eqref{chproduct} and then that
      $\partial/\partial w \contr \xi = h^{-1} \partial/\partial \eta \contr \xi$, we get
      \begin{align}\label{ch}
        \mu.\xi = 
        \int_{w=0} \sum_{\alpha \leq M} c_\alpha \partial_w^{M-\alpha}(A) \partial_w^{\alpha}(\vec \contr \xi/h).
      \end{align}
      where $M=(M_1,\dots,M_p)$ and $\vec = \partial/\partial \eta$.
      We now want to express $\phi\mu.\xi$ in terms of derivatives with respect to the variables
      $\eta_i$ instead of $w_i$. By inverting the matrix $\partial f/\partial \eta$ we get
      \begin{align}\label{rel_diff_var}
        \frac{\partial}{\partial w_j} = \frac 1h \sum_k \gamma_{jk}
        \frac{\partial}{\partial \eta_k},
      \end{align}
      where $\gamma_{jk}$ are holomorphic.
      Multiplying the test form in \eqref{ch} by $\phi$, we thus get operators $Q_\alpha$ so that
      \begin{align}\label{phimu}
        \phi\mu.\xi &= \int_{w=0} \frac{1}{h^{\sum M_j}} \sum_{\alpha \leq M}
        Q_\alpha (\phi) \partial_w^\alpha (\vec \contr \xi) = \notag\\
        &= \int_X \sum_{\alpha \leq M}\frac 1{h^{N_0}} L_\alpha (\phi) K_\alpha (\vec \contr \xi),
      \end{align}
      where $N_0$ and $N_1$ are  integers such that $N_0 = 2N_1 + \sum_jM_j$ and
      $L_\alpha = h^{N_1} Q_\alpha$ and $K_\alpha = h^{N_1} \partial^\alpha_w$ are differential operators
      with respect to the original variables $(\zeta,\eta)$. It follows from \eqref{rel_diff_var} that
      $L_\alpha$ and $K_\alpha$ are holomorphic across $W$ if $N_1$ is chosen sufficiently large.
      
      Clearly, the values of $\partial_w^\alpha (\vec \contr \xi)$ can be
      prescribed on $\{w=0\}$. Therefore $\phi \mu = 0$ on $X \setminus W$
      if and only if $L_\alpha (\phi) = 0$ on $X \setminus W$ for all $\alpha\leq M$, but by 
continuity and the $\sep$,
      these relations hold if and only if they hold across $W$. Thus $\{L_\alpha\}$ is a defining set for $\ann \mu$.

      Let now $\xi$ be an arbitrary test form.
      For the sake of simplicity,  we assume that $\phi = 1$, but the same idea works for general $\phi$.
      We thus want to show that
      \begin{align}\label{utanphi}
        \mu.\xi = \lim_{\varepsilon \to 0}\int_X \chi({|h|}^2/\varepsilon^2)
        \frac{\mathcal{L}(\vec \contr \xi)}{h^{N_0}},
      \end{align}
      where $\mathcal{L}=\sum_{\alpha \leq M} L_\alpha(1)K_\alpha$.

      The right hand side of \eqref{utanphi} defines a current $\tau$. Outside of $\{h=0\}$, $\tau$ and $\mu$
      are equal, and the latter current has the $\sep$.
      It is therefore enough to show that $\tau$ has the $\sep$.
      By expanding
      \begin{align}\label{extending}
        \chi({|h|}^2/\delta^2)\tau.\xi=\int_X \frac{\mathcal{L}(\chi({|h|}^2/\delta^2) \vec \contr \xi)}{h^{N_0}}
      \end{align}
      we get one term when all derivatives of $\mathcal{L}$ hit $\vec \contr \xi$,
      and clearly this term is precisely $\tau.\xi$ in the limit. We will now explain why all other contributions vanish;
      all these terms contain derivatives of $\chi({|h|}^2/\delta^2)$ as a factor, and
      such a factor can be written as a sum of terms
      \begin{align}\label{haakans_noll_term}
        \chi^{(k)} \left(\frac{{|h|}^2}{\delta^2}\right) \cdot \left(\frac{{|h|}^2}{\delta^2}\right)^k
        \frac{\sigma}{h^\kappa},
      \end{align}
      for some integers $k$ and $\kappa$ and a smooth function $\sigma$.
      Let us define $\tilde{\chi}(x) = \chi(x) - x^k \chi^{(k)}(x)$.
      We note that the current $\sigma/h^\kappa$ can be defined both as $\lim_{\delta \to 0} \chi({|h|}^2/\delta^2)\sigma/h^\kappa$ and
      $\lim_{\delta \to 0} \tilde{\chi}({|h|}^2/\delta^2)\sigma/h^\kappa$, so it follows that \eqref{haakans_noll_term} must be zero
      in the limit.
    \end{proof}
  \end{theorem}




  
  \section{Residue currents associated to ideals}\label{currents}
  We will give a summary of the machinery of residue currents needed to prove Theorem~\ref{bs}.
  In \cite{anderssonwulcan:prescribed} was presented a method for constructing currents $R$ and $U$
  associated to any generically exact complex
  \begin{align}\label{general_complex}
    \dots \to E_2 \overset{f_2}{\to} E_1 \overset{f_1}{\to} E_0 
  \end{align}
  of hermitian vector bundles over an open set in $\C^n$. The total bundle is then $(E,f)$, where
  $E=\bigoplus E_k$, $f=\oplus f_j$.
  The currents $U$ and $R$ take values in the bundle $\End E$. Let $\sigma_k: E_{k-1} \to E_k$
  be the mapping of pointwise minimal norm that is the inverse of $f_k$ on the image of $f_k$ and extend it by zero on the orthogonal
  complement of the image. Let $W$ be the analytic set where \eqref{general_complex} is not exact. Outside of $W$, we define
  \begin{align}\label{u_def}
    u &= \sum_{j=1}^{n+1} u_j \\\notag
    u_j &=  \dbar \sigma_j \wedge \ldots \wedge (\dbar \sigma_2) \sigma_1,
  \end{align}
  so that $u_j$ is the $(0,j-1)$-bidegree component of $u$ which takes values in $\Hom(E_0,E_j)$.
  One can then extend $u$ to a current $U$ across $W$ by setting
  \begin{align}\label{U_def}
    U = \lim_{\epsilon \to 0} \chi(|F|^2/\epsilon^2)\wedge u,
  \end{align}
  where $F$ is a holomorphic tuple vanishing on $W$. 
  We define an operator $\nabla_f$ acting on currents with values in $E$ by
  \begin{align}\label{nabla_def}
    \nabla_f &= f - \dbar.
  \end{align}
  The residue current $R$ is defined by
  \begin{align}\label{UR_relation}
    \nabla_{f} \circ U = 1 - R.
  \end{align}  
  One can check that $\nabla_{f} \circ U = 1$ on $X\setminus W$, so $R$ has support on $W$.
  The $(0,k)$-bidegree component of $R$ is denoted by $R_k$. 
  Since $\nabla_{f}^2 = 0$, \eqref{UR_relation} gives that $\nabla_{f} R = 0$, which is equivalent to
  \begin{align}\label{nabla_written_out}
    f_1 R_1 &= 0 \\\notag
    f_{k+1} R_{k+1} &= \dbar R_k, \quad k\geq 1.
  \end{align}

  An easy calculation shows that
  \begin{align}\label{R_reg}
    R=\lim_{\epsilon \to 0} R_{0,\epsilon} + R_{1,\epsilon} + \ldots + R_{n,\epsilon},
  \end{align}
  where
  \begin{align}\label{R_reg_comp1}
     R_{0,\epsilon}&=(1-\chi(|F|^2/\epsilon^2))\\\notag
     R_{j,\epsilon}&= \dbar\chi(|F|^2/\epsilon^2) \wedge u_j, \quad 1\leq j \leq n.
  \end{align}  

  We restrict our attention to the case when $\rank E_0 = 1$.
  For any ideal $\idealJ \subset \local$ we can choose \eqref{general_complex} (in various ways) so that
  $\idealJ = \Image(\O(E_1) \to \O(E_0))$.
  Then $\ann R \subset \idealJ$; indeed, if $\phi R = 0$, then by \eqref{UR_relation} $\nabla_f U \phi = \phi$. By solving a sequence
  of $\dbar$-equations one can then show that $\phi \in \idealJ$, see \cite{andersson:residueideals}. The converse inclusion does not
  hold in general.
  A main result of \cite{anderssonwulcan:prescribed}
  is that if \eqref{general_complex} is chosen so that $\O(E_k)$ with maps $(f_k)$ is a resolution of $\O_{\C^n}/\idealJ$, that is,
  an exact complex of sheaves, then $\ann R = \idealJ$.
  Another choice for \eqref{general_complex} is the Koszul complex. This has the advantage that the residue has an explicit form, and
  if $\idealJ$ happens to be a complete intersection, then the Koszul complex is a resolution, so $\ann R = \idealJ$.

  Now let $\idealJ \subset\local$ be the ideal that defines the analytic space $Z$ of Theorem~\ref{bs}.
  Thus we assume that $\idealJ$ has pure codimension $p$. We can then
  choose \eqref{general_complex} so that the corresponding sheaf complex is a resolution of $\local / \idealJ$.
  Let $R^Z$ be the current associated to $\idealJ$ as above so that $\ann R^Z = \idealJ$.
  It has support on $W=X$. 

  We also need a current $R^{a^r}$ with the property $\ann R^{a^r} \subset (a)^r$. We obtain such a current by the procedure
  above, where the complex \eqref{general_complex} should be chosen so that $(a)^r = \Image(\O(E_1) \to \O(E_0))$.
  The form $u$ and the current $U$ will be denoted by $u^{a^r}$ and $U^{a^r}$ respectively.
  For the actual choice of \eqref{general_complex},
  we follow \cite{andersson:bsexplicit}, where $u^{a^r}$ and $R^{a^r}$ are described explicitly.
  Define
  \begin{align}\label{sigma_def}
    \sigma_i = \sum_{j=1}^m\frac{\overline{a_j}e^i_j}{|a|^2}
  \end{align}
  outside of $W=Z(a)$, where $\{e^i_j\}_j$ are frames of some trivial vector bundles.
  Andersson shows, see eq.~(2.3) and the beginning of Section~3 in \cite{andersson:bsexplicit}, that
  \begin{align}\label{u_a_r}
    u^{a^r} = \sum_{\substack{\sum j_i \leq \varrho -1}}  \sigma^1\wedge
    \left(\dbar \sigma^1 \right)^{\wedge j_1}\wedge\dots\wedge
    \sigma^r\wedge \left(\dbar \sigma^r \right)^{\wedge j_r},
  \end{align}
  where $\varrho = \min(n-p,m)$.
  We then obtain $R^{a^r}$ by \eqref{R_reg} and \eqref{R_reg_comp1} with $F = |a|$, that is,
  \begin{align}\label{R_a_r}
    R^{a^r} = \lim_{\epsilon \to 0} (1-\chi^a_\epsilon) + \dbar \chi^a_\epsilon \wedge u^{a^r},
  \end{align}
  where $\chi^a_\epsilon = \chi(|a|^2/\epsilon^2)$.

  
  \subsection{Almost semi-meromorphic currents}\label{nsm}

  In order to describe the structure of the current $R^Z$ we need to make a digression
  into {\it almost semi-meromorphic} currents.
  Recall that a current is semi-meromorphic if is a principal value
  current of the
  form $\alpha/f$, where $f$ is a holomorphic section of a line bundle and $\alpha$ is a smooth 
  section of the same bundle.
  We say that a current $a$ in an open subset $V\subset \C^n$ is almost semi-meromorphic
  if there is a modification $\pi\colon V'\to V$ and                
  a semi-meromorphic current $\tilde a$ on $V'$ such that
  $a=\pi_*\tilde a$.

  Notice that if  $a$ is almost semi-meromorphic and $\xi$ is smooth,
  then $\xi\wedge a$ is almost semi-meromorphic. In fact,
  $\xi\w a=\pi_*(\pi^*\xi\w\tilde a)$.


  We will have use for the following result which is a part of Theorem~5.1 in \cite{AW3}.            

  \begin{lemma}\label{nsm1} If $a$ is an almost semi-meromorphic 
    $(0,*)$-current in an open subset of $\C^n$, 
    then also $(\partial/\partial z_\ell) a$ is almost semi-meromorphic.
  \end{lemma} 




By \cite[Theorem~4.6]{AW3} one can define the product $a\w \mu$,
where $a$ is almost semi-meromorphic and $\mu$ is either an 
almost semi-meromorphic current or a Coleff-Herrera current, in the following way; notice first that $a$ is smooth outside of
an exceptional analytic set. If $\mu$ is almost semi-meromorphic, let $\{h=0\}$ be any hypersurface
containing the exceptional set of $a$ and let $a\w \mu:=\lim_{\epsilon\to 0}\chi(|h|^2/\epsilon^2)a\w \mu$;
the resulting current is independent of the choices of $\chi$ and $h$ and it is 
almost semi-meromorphic. If $\mu$ is a Coleff-Herrera current on the analytic set $W$
we will assume that the exceptional set of $a$ intersects $W$ properly since this is the only case of interest for us.
Then let $\{h=0\}$ be a hypersurface containing the exceptional set of $a$ such that $h$ is generically non-vanishing on $W$
and as before set $a\w \mu:=\lim_{\epsilon\to 0}\chi(|h|^2/\epsilon^2)a\w \mu$. Again the resulting current is independent
of the choices made and, by \cite[Corollary~4.7]{AW3}, $a\w \mu$ is supported on $W$ and has the SEP with respect to $W$.


  We also need a result which is a slight variation of (part of) Proposition~3.3 in \cite{as_grothendieck}.
  We include here the needed part of the proof.
  \begin{proposition}\label{RZfactored}
    There is a vector-valued Coleff-Herrera current $\mu$ and an
    endomorphism-valued almost semi-meromorphic current $b$ in a neighborhood of $0$ in $\C^n$ such that
    $b$ is smooth outside of $X_{sing}$ and such that $R^Z = b \mu$.  
  \end{proposition}  
    \begin{proof}
      Decompose $R^Z = \sum_{k=0}^n R^Z_k$ so that $R^Z_k$ has bidegree $(0,k)$. Proposition~2.2 in \cite{anderssonwulcan:prescribed}
      gives that the conjugate of $I(Z) = \sqrt{\idealJ}$ annihilates $R^Z$ and that $R^Z_k = 0$ for $k<p$, so 
      $R^Z=R^Z_p + R^Z_{p+1} + \ldots$. 
      Consider the dual complex of \eqref{general_complex}
      \begin{align*}
        0 \to \O(E^*_0) \overset{f^*_1}{\to} \O(E^*_1) \to \dots \to \O(E^*_p) \overset{f^*_{p+1}}{\to} \O(E^*_{p+1}) \to \dots,
      \end{align*}
      where $f^*_j$ is the matrix transpose of $f_j$; we now assume that the sheaf complex associated with \eqref{general_complex}
      is a resolution of $\O/\mathcal{J}$. Since the sheaf $\kersh f^*_{p+1}$ is coherent, we have an exact
      sequence
      \begin{align*}
        \O(F^*) \overset{g^*}{\to} \O(E^*_p) \overset{f^*_{p+1}}{\to}\O(E^*_{p+1}),
      \end{align*}
      where $F$ is a vector bundle. Since $\nabla_f R^Z = 0$, we have $\dbar R^Z_p = f_{p+1} R^Z_{p+1}$, which gives that
      \begin{align*}
        \dbar (g R^Z_p) = g \dbar R^Z_p = g f_{p+1} R^Z_{p+1} = 0,
      \end{align*}
      since $gf_{p+1}=0$. Since, by Corollary~2.4 in \cite{anderssonwulcan:decomp} and Proposition~2.2 in 
\cite{anderssonwulcan:prescribed}, $R^Z_p$ has the SEP 
with respect to $X$ and is annihilated by anti-holomorphic functions vanishing on $X$,
it follows that $gR^Z_p$ is a Coleff-Herrera current on $X$.
      
      Let $Z_{p+1}$ be the complement of the set where $f_{p+1}$ has optimal rank.
      This is an analytic subset of $X_{sing}$ that is intrinsic to $Z$. Outside of $Z_{p+1}$, the mapping $g: E_p \to F$ has
      constant rank. We define a mapping $\sigma_F: F \to E_p$ on the complement of $Z_{p+1}$ by
      \begin{align*}
        (\sigma_F)_{|(\Image g)^{\bot}} &= 0 \\
        \sigma_{F}g_{|(\ker g)^{\bot}} &= 1_{(\ker g)^{\bot}}.
      \end{align*}
      This means that $\sigma_F$ is the minimal norm inverse of $g$ on the image of $g$, and it is
      extended by zero to the orthogonal complement of this image; we here choose an auxiliary Hermitian metric on $F$.
      It is shown in Section~2 of \cite{anderssonwulcan:prescribed} that $\sigma_F$ has an almost semi-meromorphic extension
      across $Z_{p+1}$; the extension is denoted $\sigma_F$ as well.
      From equation \eqref{u_def}, we see 
      that outside of $Z_{p+1}$, $u^Z_p$ takes values in the subbundle $(\Image f_{p+1})^{\bot}=(\ker g)^{\bot} \subset E_p$ since
      this is true for $\sigma^Z_p$. It follows that $R^Z_p = \sigma_F g R^Z_p$ outside of $Z_{p+1}\subset X_{sing}$
and since both $R^Z_p$ and $\sigma_F g R^Z_p$ have the SEP with respect to $X$ we see that $R^Z_p = \sigma_F g R^Z_p$
holds cross $Z_{p+1}$.

      
      Now, let $Z_{p+l}$ be the set where $f_{p+l}$ does not have optimal rank; this is again an 
      analytic subset of $X_{sing}$ that is intrinsic to $Z$ and since $Z$ has pure dimension
      we also have that $\textrm{codim}\, Z_{p+l} \geq p+l+1$, see Corollary 
~20.14 in \cite{Eis}. By Theorem~4.4 in \cite{anderssonwulcan:prescribed} there are almost
      semi-meromorphic endomorphism valued forms $\alpha_l$ such that $R^Z_{p+l} = \alpha_l R^Z_p$ outside of
      $Z_{p+l}$. As in the proof of Proposition~3.3 in \cite{as_grothendieck} this equality extends across $Z_{p+l}$
      and we are done.
    \end{proof}
  

  \section{Proof of Theorem~\ref{bs}}\label{proof}

  We will use the product of the currents $R^{a^r}$ and $R^Z$ which was defined in \cite{ass:bsvar}: 
  \begin{align}\label{Rar_RZ_prod}
    R^{a^r} \wedge R^Z = \lim_{\epsilon \to 0} \left[R_{0,\epsilon}^{a^r} + \dbar\chi^a_\epsilon \wedge u^{a^r}\right] \wedge R^Z.
  \end{align}
  This current takes values in the tensor product of the two corresponding complexes \eqref{general_complex} for $R^{a^r}$ and $R^Z$.

  It follows by Proposition 2.2 in \cite{ass:bsvar} (which holds in the non-reduced setting with the same proof),
  that $\phi R^{a^r} \wedge R^Z = 0$ for $\phi \in \local$
  implies that $\phi \in \idealJ + (a)^r$, that is, the image of $\phi$ in $\localz$ belongs to $\ideala^r$.
  Although $R^{a^r} \wedge R^Z$ is a current in $\C^n$, $\phi R^{a^r} \wedge R^Z$ depends only on the image 
  of $\phi$ in $\local /\idealJ$; in fact,
  if $\phi \in \idealJ$, then $\phi R^Z =0$ and so 
  $\phi R^{a^r} \wedge R^Z = \lim_{\epsilon\to 0} \phi R^{a^r}_\epsilon \wedge R^Z = 0$. 
  Moreover, $R^{a^r} \wedge R^Z$ only depends on the images of the generators $a_j$ in $\local / \idealJ$.
  This can be deduced for example
  from Proposition~2.2 in \cite{anderssonwulcan:prescribed}.
  Hence the proof of Theorem \ref{bs} is reduced to
  showing that $\phi R^{a^r} \wedge R^Z = 0$ if $\phi$ satisfies \eqref{premise} with a suitable constant $N$.

  By Proposition~\ref{RZfactored}, each component of $R^Z$ is a sum of terms that can be factored
  as $b \mu$ where $b$ is almost semi-meromorphic
  and $\mu$ is a Coleff-Herrera current. According to Theorem~\ref{bjork_thm}, the annihilator of each such $\mu$ has
  a defining set $M^\mu = \{L_1^\mu, \dots, L_{M_\mu}^\mu\}$ (and these operators satisfy \eqref{phimu'}).
  As the operators $L_1, \dots, L_M$ in Theorem~\ref{bs},
  we take the union of $M^\mu$ over all Coleff-Herrera currents $\mu$ that arise in the way we described.
  
  Although
  \begin{align}\label{to_be_zero_all_terms}
    \phi R^{a^r} \wedge R^Z=\lim_{\epsilon \to 0} \phi (R^{a^r}_{0,\epsilon} + \dbar \chi^a_\epsilon \wedge u^{a^r}) \wedge R^Z,
  \end{align}
  we will only prove that $\lim_{\epsilon \to 0} \phi \dbar \chi^a_\epsilon \wedge u^{a^r} \wedge R^Z =0$ as the proof for the remaining term
  is similar but easier.
  It suffices to show that
  \begin{align}\label{want_to_show}
    \phi \dbar \chi^a_\epsilon \wedge u^{a^r}\wedge b\mu.\omega
  \end{align}
  has limit $0$ as $\epsilon$ goes to $0$, where $\omega$ is a test form. 
  Let $h$ be the holomorphic function defined in \eqref{detdef} and set 
  $\chi^b_\delta = \chi(|h|^2/\delta^2)$. 
  We apply \eqref{phimu'} with 
  $\xi=\chi^b_\delta b\, \dbar \chi^a_\epsilon \wedge u^{a^r} \wedge \omega$ to compute \eqref{want_to_show}.
  This gives
  \begin{align}\notag
    \phi \dbar \chi^a_\epsilon \wedge u^{a^r}\wedge b\mu.\omega &=
    \lim_{\delta\to 0} \phi \mu.\dbar \chi^a_\epsilon \wedge u^{a^r} \wedge \chi^b_\delta b\omega =\\\label{phimuRa}
    &= \lim_{\delta\to 0} \int_X \sum_{\alpha \leq M}
    \frac {1}{h^{N_0}} L_\alpha (\phi) K_\alpha (\vec \contr \dbar \chi^a_\epsilon \wedge u^{a^r} \wedge \chi^b_\delta b\omega).
  \end{align}
  
  \begin{proposition}\label{nsm_op_lma}
    There are differential operators $\tilde{K}_\alpha$ with almost semi-meromorphic
    coefficients and holomorphic derivatives so that
    \begin{align}\label{nsm_operators}
      \phi \dbar \chi^a_\epsilon \wedge u^{a^r}\wedge b\mu.\omega = 
      \int_X \sum_{\alpha \leq M} L_\alpha (\phi) \tilde{K}_\alpha (\dbar \chi^a_\epsilon \wedge u^{a^r}\wedge \omega).
    \end{align}
  \end{proposition}

  This integral should be interpreted as a principal value, that is, as the limit of
  \begin{align*}
    \int_X \sum_{\alpha \leq M} L_\alpha (\phi) \chi^b_\delta \tilde{K}_\alpha (\dbar \chi^a_\epsilon \wedge u^{a^r}\wedge \omega)
  \end{align*}
  as $\delta \to 0$.

  \begin{proof}
    We apply Leibniz' rule to see that the integral in \eqref{phimuRa} is a sum of terms of the form  
    \begin{align}\label{kut}
    \int_X \psi\wedge \partial^{\alpha}\big(\dbar \chi_{\epsilon}^a\wedge u^{a^r}\wedge \omega\big) \wedge
    \partial^{\beta}\chi^b_{\delta}\wedge \frac{\partial^{\gamma} b}{h^{N_0}},
    \end{align}
    where $\psi$ is smooth. By Lemma~\ref{nsm1} and the discussion  following it we have that $\partial^{\gamma} b/h^{N_0}$
    is almost semi-meromorphic. Moreover, it is smooth outside $\{h=0\}$ and
    since $[X]$ is a Coleff-Herrera current on $X$ (with values in the $(p,0)$-forms) it follows that
    \begin{align*}
    \frac{\partial^{\gamma} b}{h^{N_0}}\wedge [X] = \lim_{\delta\to 0} \chi_{\delta}^b \frac{\partial^{\gamma} b}{h^{N_0}}\wedge [X].
    \end{align*}
    Hence, if $\beta=0$ then \eqref{kut} goes to an expression of the form \eqref{nsm_operators} as $\delta\to 0$. 
    Assume now that $\beta\neq 0$. Then $\partial^{\beta}\chi_{\delta}^b$ is a sum of terms of the form \eqref{haakans_noll_term}
    with $k\geq 1$. It follows in the same way as in (the end of) Section~2 that, in this case, \eqref{kut} goes to $0$ as $\delta\to 0$.

  \end{proof}

  We now choose a resolution $X' \overset{\pi}{\to} X$ such that $X'$ is smooth and $\pb h$ is locally a monomial
  and the coefficients of $\tilde{K}_\alpha$ are push forwards of semi-meromorphic forms whose denominators have
  normal crossing zero sets.
  It suffices to show that
  \begin{align}\label{suffices}
    I_\epsilon = \int_{X'} \frac{\xi'}{s_1^{1+n_1}\cdot\ldots\cdot s_{n-p}^{1+n_{n-p}}}\wedge
    \pb{\left(L_\alpha \phi\right)}\pb\left(\partial^\beta_{\eta}(\dbar \chi^a_\epsilon \wedge u^{a^r})\right) \to 0,
  \end{align}
  where $\beta \leq M$, $\xi'$ is a smooth form with compact support and $\eta$ is a local coordinate system on the ambient space $\C^n$ of $X$.
  This integral is defined as a principal value as before.
  
  We want to integrate by parts in \eqref{suffices}.
  The principal value current
  \begin{align*}
    \left[1/{s_1^{1+n_1}\cdot\ldots\cdot s_{n-p}^{1+n_{n-p}}} \right]
  \end{align*}
  is a tensor product of one variable currents
  $\left[1/{s_j^{1+n_j}} \right]$.
  Furthermore, one has for $m \geq 1$,
  \begin{align*}
    \frac{\partial}{\partial s_j}\left[\frac{1}{s_j^m}\right] = -m \left[\frac{1}{s_j^{m+1}}\right].
  \end{align*}
  This yields indeed that
  \begin{align}\label{ideltaparts}
    I_\epsilon = \int_{X'} \frac{ds}{s}\wedge
    \partial_s^{(n_1,\dots,n_{n-p})}\left(\xi'\wedge\pb{\left(L_\alpha \phi\right)}\pb\left(\partial^\beta_{\eta}(
    \dbar \chi^a_\epsilon \wedge u^{a^r})\right)\right),
  \end{align}
  where
  \begin{align*}
    \frac{ds}{s} = \frac{ds_1 \wedge\dots\wedge ds_{n-p}}{s_1\cdot\dots\cdot s_{n-p}}.
  \end{align*}

  As the constant $N$ in the formulation of Theorem~\ref{bs} we set
  \begin{align}\label{sufficient_N}
    N=\max \left[2\varrho + |M| + \sum_{j=1}^{n-p} n_j \right], 
  \end{align}
  where $\varrho = \min(m,\dim X) = \min(m,n-p)$ and the quantities $|M|$ and $\sum_1^{n-p} n_j$ are
  maximized over the components of $R^Z$ and over all local charts of a (finite) covering of $X'$ on which $h$ is a monomial.
  We claim that if \eqref{premise} holds, then for any integers $k_j \leq n_j$, $1\leq j\leq n-p$,
  \begin{align}\label{remains}
    \partial_s^{(k_1,\dots,k_{n-p})}\left(\pb{\left(L_\alpha \phi\right)}\pb\left(\partial^\beta_{\eta}(
    \dbar \chi^a_\epsilon \wedge u^{a^r})\right)\right)
  \end{align}
  is bounded by a constant that is independent of $\epsilon$. Note that the pointwise limit of any derivative of $\dbar \chi^a_\epsilon$
  is zero almost everywhere.
  An application of dominated convergence in \eqref{ideltaparts} thus gives that $I_\epsilon \to 0$,
  and thereby concludes the proof, given the claim.

  It remains to prove the claim above. We will assume that $k_j=n_j$ for $1\leq j \leq n-p$, as this is the worst case.
  We extend $\pi$ to a resolution $X'' \overset{\pi'}{\to} X' \overset{\pi}{\to} X$ 
  that principalizes $(a)$, that is
  $\pbpi\pb a_j = a_0 a'_j$ for some non-vanishing tuple $\{a'_j\}$.
  Note that the form
  $ds/s$ becomes a sum of similar forms when we pull it back to $X''$.
  After pulling back \eqref{remains} to $X''$, we write it as a linear combination of forms
  \begin{align}\label{remains2}
     \pbpi \left[\partial_s^{(\hat{n}_1,\dots,\hat{n}_{n-p})} \pb{(L_\alpha \phi)} \right]
     \pbpi \left[\partial_s^{(\tilde{n}_1,\dots,\tilde{n}_{n-p})} \pb\left(\partial^\beta_{\eta}(
      \dbar \chi^a_\epsilon \wedge u^{a^r})\right)\right],
  \end{align}
  for integers $\hat{n}_j$ and $\tilde{n}_j$ such that $\hat{n}_j + \tilde{n}_j = n_j$.
  
  By \eqref{premise} we have that
  \begin{align}\label{sz_trick_prem}
    |\pb{(L_\alpha \phi)}| \lesssim |\pb a|^{N+r-1}.
  \end{align}
  The classical Brian\c con-Skoda theorem therefore implies that locally $\pb{(L_\alpha \phi)} \in~(\pb a)^{N+r-\varrho}$.
  For the first factor of \eqref{remains2}, we therefore locally have that
  \begin{align}\label{L_phi_der}
    \pbpi \left[\partial_s^{(\hat{n}_1,\dots,\hat{n}_{n-p})} \pb{(L_\alpha \phi)} \right] \in (a_0)^{N+r-\varrho -\hat{n}},
  \end{align}
  where $\hat{n} = \hat{n}_1 + \ldots + \hat{n}_{n-p}$.
  As for the second factor of \eqref{remains2}, we argue that
  \begin{align}\label{res_der}
    \pbpi \left[\partial_s^{(\tilde{n}_1,\dots,\tilde{n}_{n-p})} \pb\left(\partial^\beta_{\eta} (\dbar \chi^a_\epsilon \wedge u^{a^r})\right)\right] =
    \O(|a_0|^{-(|\beta| + \tilde{n} + \varrho + r)}),
  \end{align}
  where $\tilde{n} = \tilde{n}_1 + \ldots + \tilde{n}_{n-p}$. The chain rule gives that
  \begin{align*}
    \partial_{s_j} \pb \omega = \sum_{j=1}^n (\partial_{s_j}\pi_j) \pb \partial_{\eta_j} \omega,
  \end{align*}
  for any $(0,q)$-form $\omega$. Consequently, we can rewrite the left hand side of \eqref{res_der} as a sum of terms like
  \begin{align}\label{from_res_der_simpler}
    \xi'' \pbpi \pb \left[\partial^{\tilde{\beta}}_{\eta} (\dbar \chi^a_\epsilon \wedge u^{a^r})\right],
  \end{align}
  where $\xi''$ is a smooth function and $|\tilde{\beta}| = |\beta| + \tilde{n}$. A slight reformulation of Lemma~4.2
  in \cite{sznajdman:ar} states that
  \begin{lemma}\label{reformulation}
    Let $\tilde{\pi}: X'' \to X$ be any principalization of the ideal $(a)$ on $X$ so that
    $\tilde{\pi}^* (a) = (a_0)$. Then for any multi-index $\beta_0$
    \begin{align*}
      \tilde{\pi}^* \left[\partial^{\beta_0}_{\eta} (\dbar \chi^a_\epsilon \wedge u^{a^r})\right] =
      \sum_{j=0,1} \left(\frac{d\overline{a_0}}{\overline{a_0}}\right)^j \wedge \O(|a_0|^{-(|\beta_0| + \varrho + r-1)}).
    \end{align*}
  \end{lemma}
  We apply this lemma to \eqref{from_res_der_simpler} with $\tilde{\pi} = \pi \circ \pi'$ and $\alpha_0=\tilde{\beta}$. We
  will use that $d\overline{a_0}/\overline{a_0}$ is integrable in the following section, but here we can just estimate
  $d \overline{a_0}$ by $\O(1)$.
  We conclude by \eqref{L_phi_der} and \eqref{res_der} that the form in \eqref{remains} is indeed bounded if
  \begin{align}\label{N_condition}
    N \geq |\beta| + 2 \varrho + \hat{n} + \tilde{n} = |\beta| + 2\varrho + \sum_{j=1}^{n-p} n_j,
  \end{align}
  and $|\beta| \leq |M|$.

\begin{remark}\label{global} Let $Z\subset{\mathbb P}^n$ be a pure-dimensional projective variety,
let $\idealJ$ be the associated coherent ideal sheaf and let $X$ be the
associated reduced projective variety $X$. We believe that one can
globalize Proposition~2.2 above and obtain
Noetherian operators $L_{\alpha}$ in $\C^n$ such that an analogue of Theorem 1.2 holds: There is a constant
$N$ such that if $\ideala$ is a polynomial ideal and (1) holds (with say $r=1$), then $\phi$ belongs to the
ideal $\ideala$.  

In general, the fact that $\phi$ is in $\ideala$ only implies  that there is a representation
$\phi=\sum a_jq_j$, where $\deg a_j q_j$ is like $\deg \phi+2^{2^d}$, if 
$d$ is the (maximal) degree of the generators of $\ideala$.  However, in the reduced case,
i.e., as in \cite{aw:semester},  the condition (1) implies a  degree estimate
like $\deg\phi+Cd^{n-p}$. One could hope for a result of this kind even in the
global non-reduced case. 
\end{remark}

\begin{remark}\label{B-Sato}
Instead of using a log resolution $X'\to X$ and integrate by parts on $X'$ to see that
the integral in \eqref{nsm_operators} is a sum of terms of the form \eqref{ideltaparts}
it would be interesting to try to
proceed as follows: Make successive normalized blow-ups along suitable ideals
to get a normal modification $X'\to X$ such that the almost semi-meromorphic 
coefficients of all the operators $\tilde{K}_{\alpha}$ are push-forwards of semi-meromorphic
currents on $X'$. Then use Bernstein-Sato functional equations on $X'$ to see that
the integral in \eqref{nsm_operators} is a sum of terms of a form very similar to \eqref{ideltaparts}.
In cases where the orders of the differential operators appearing in the Bernstein-Sato equations
are known one could thus possibly
be able to estimate the number $N$ in Theorem~\ref{bs}.
\end{remark}

  \section{The Cohen-Macaulay case}\label{smooth}
In general it is hard to get a good estimate of the constant $N$ in Theorem~\ref{bs} from  the
proof. In the case when $Z$ is Cohen-Macaulay things are somewhat simpler.
To begin with, then  $R^Z$ only consists of the term $R^Z_p$ which is a 
vector-valued Coleff-Herrera current on $X$, see \cite{andersson:chduality}.
That is, the factor $b$ in
Proposition~3.3 is just $1$ and so one does not need the technical Proposition~\ref{nsm_op_lma}
to prove Theorem~\ref{bs} in the Cohen-Macaulay case; it is sufficient to use Theorem~\ref{bjork_thm}.

To have the representation (11), with $\mu=R^Z_p$,
it is precisely 
required that $f_j^{1+M_j}$ are in $\idealJ$,  and thus $M_j$ are governed
by the constant in the local Nullstellensatz. 
If, for instance, the generators are polynomials,
then one can control $M_j$ by the degree of the polynomials via Bezout estimates.
Notice  that $M_j$ give an upper
bound of the degree of the resulting Noetherian operators $L_j$.
It could be mentioned here  that
in the Cohen-Macaulay case also the function $A$ in (11) can be obtained more 
explicitly  by means of a formula from \cite{lark}. The function $h$ in the denominator
in (8) is related to the structure form in \cite{as_grothendieck},
and the power $N_0$ depends on the degree of the resulting Noetherian operators. 

\smallskip

  We consider now the case where in addition $X$ is smooth.
Choose the smallest possible multi-index $M=(M_1,\ldots, M_p)$ as above. 
For any component $\mu_j$, there is a defining set of operators $\{L_{j,\alpha}\}$, $\alpha\le M$, such that
\eqref{phimu'} holds. The union  $\{L_{j,\alpha}\}$ of  these sets is a defining set for $\idealJ$. 
Let us introduce the  notation $L_\alpha$ for the vector-valued operator $L_{j,\alpha}$.
 
 Notice that, since $X$ is smooth,
  we can take $w=\eta$ in \eqref{coordinates} so that $h=1$ in \eqref{detdef}.

  For each differential operator $L$ we denote by $\order(L)$ its order as a differential operator in the ambient space.
  One can also define the order of an element in $\dmod$ as the minimal order of any operator representing it.

  Let $d$ be the maximal distribution order all Coleff-Herrera currents
  that are annihilated by $\idealJ$. One can show that this number does not depend on the embedding
  of $Z$ into $\C^n$. By Example 1 in \cite{andersson:chduality},
  the components of $\mu=R^Z$ generate the Coleff-Herrera currents annihilated by $\idealJ$.
  Thus $d$ can also be expressed as the maximal distribution order of the components of $R^Z$.

  \begin{theorem}\label{bssmooth}
    Let $Z$ be a germ of a Cohen-Macaulay analytic space such that $X$ is smooth,
    and let $\{L_\alpha\}$, $\alpha\le M$, be the vector-valued Noetherian operators obtained from
    $\mu=R^Z_p$ as above.
    If $r\geq 1$, and $\ideala \subset \localz$ can be generated by $m$ elements, then
    \begin{align} \label{premisesmooth}
      |L_\alpha \phi|\leq C |\ideala|^{\min(m,\dim X)+d-\order(L_\alpha)+r-1} \,\,\, \textrm{on} \,\,\, X, \quad \alpha \leq M,
    \end{align}
    implies that $\phi\in\ideala^r$.
  \end{theorem}
  For the analytic space in Example~\ref{bevisatexempel}, \eqref{premisesmooth} coincides with the optimal
  hypotheses \eqref{smpremise}.
                                              
  \begin{proof}
    We can assume that $X = \C^{n-p} \subset \C^n$, and we call the last $p$ coordinates
    $w_1, \dots, w_p$. Clearly these functions form a complete intersection, so if $\mu=R^Z_p$, we have
    \begin{align*}
      \mu = A\left[\dbar \frac 1{w_1^{1+M_1}} \wedge\dots\wedge
	\dbar \frac 1{w_p^{1+M_p}} \right],
    \end{align*}
    for some (vector-valued) holomorphic function $A$, cf.\ \eqref{chproduct}. A basic computation 
    rule for Coleff-Herrera products is that
    \begin{align}\label{ch_rule}
      w^\alpha \left[\dbar \frac 1{w^{1+M}} \right] =  \left[\dbar \frac 1{w^{1+M - \alpha}} \right],
    \end{align}
    where we have used multi-index notation. Using this, we get
    \begin{align}\label{mu_expansion}
      \mu = \sum_{\alpha \leq M} C_\alpha A_\alpha(z) \left[ \dbar \frac 1{w^{1+\alpha}}\right],
    \end{align}
    where $z=(z_1,\dots,z_{n-p})$, and $A_\alpha$ are holomorphic functions and $C_\alpha$ suitable constants
    so that
    \begin{align}\label{glatt_mu_rep}
      \mu.\xi = \int_{w=0} \sum_{\alpha \leq M}
      A_\alpha(z) \partial_w^{\alpha}(\partial/\partial w \contr \xi).
    \end{align}
Notice that since $M$ is minimal, $A_M$ must be nonzero.  Therefore,
the distribution order of $\mu$ is $|M|$.
    
   
    As in Section~\ref{coleffherrera}, we multiply \eqref{glatt_mu_rep} by $\phi$ to obtain
    \begin{align}\label{smoothphimu}
      \phi\mu.\xi &= \int_{w=0} \sum_{\alpha \leq M} 
      L_\alpha (\phi) \partial_w^\alpha (\partial/\partial w \contr \xi),
    \end{align}
    where
    \begin{align*}
      L_\alpha = \sum_{M \geq \gamma\geq \alpha} {{\gamma} \choose {\alpha}}
      A_\gamma(z) \partial_w^{\gamma-\alpha}. 
    \end{align*}
    Note that
    \begin{align}\label{L_alpha_order_bound}
      \order(L_\alpha) \leq |M|-|\alpha|.
    \end{align}
    We substitue $\xi=\dbar \chi^a_\epsilon \wedge u^{a^r} \wedge \omega$ into \eqref{smoothphimu}.
    The analogue of \eqref{suffices} now becomes
    \begin{align}\label{ideltaparts_smooth}
      I_\epsilon = \int_{w=0}  L_\alpha (\phi) \partial^\alpha_w (\dbar \chi^a_\epsilon \wedge u^{a^r}) \wedge \xi'.
    \end{align}

    Applying again Lemma~\ref{reformulation}, we get
    \begin{align}\label{res_der_smooth}
      \pbpi \partial^\alpha_w (\dbar \chi^a_\epsilon \wedge u^{a^r}) =
      \sum_{j=0,1} \left(\frac{d\overline{a_0}}{\overline{a_0}}\right)^j \wedge \O(|a_0|^{-(|\alpha| + \varrho + r - 1)}).
    \end{align}
    By \eqref{L_alpha_order_bound}, we have $d - \order(L_\alpha) \geq |M| - \order(L_\alpha) \geq |\alpha|$.
    Thus, assuming that $|L_\alpha(\phi)|\leq C|a|^{\varrho+d-\order(L_\alpha)+r-1}$, it follows from \eqref{res_der_smooth} and
    the fact that $d\overline{ a_0}/\overline{a_0}$ is integrable that
    dominated convergence may be applied in \eqref{ideltaparts_smooth}. Since $\dbar \chi^a_\epsilon$ and all of its derivatives
    go to zero almost everywhere, we see that
    $\lim_{\epsilon \to 0} I_\epsilon = 0$, which was to be shown.
  \end{proof}
    
{\bf Acknowledgements}: I thank my advisors Mats Andersson and H\aa kan Samuelsson who gave insightful comments on the many preliminary versions.

  \providecommand{\bysame}{\leavevmode\hbox to3em{\hrulefill}\thinspace}
\providecommand{\MR}{\relax\ifhmode\unskip\space\fi MR }
\providecommand{\MRhref}[2]{%
  \href{http://www.ams.org/mathscinet-getitem?mr=#1}{#2}
}
\providecommand{\href}[2]{#2}

\end{document}